\documentclass[10pt]{amsart}
\usepackage{amssymb, amsmath, amsthm, enumerate, times}

\newtheorem*{staircase}{Generalized Climbing Stairs Problem}
\newtheorem*{partition}{Theorem}
\newtheorem*{order}{Generalized Climbing Stairs Problem with Order}

\newtheorem{question1}{Question}
\newtheorem{question2}[question1]{Question}
\newtheorem{question3}[question1]{Question}
\newtheorem{question4}[question1]{Question}
\newtheorem{question5}[question1]{Question}
\newtheorem{question6}[question1]{Question}
\newtheorem{question7}[question1]{Question}
\newtheorem{question8}[question1]{Question}
\newtheorem{question9}[question1]{Question}
\newtheorem{question10}[question1]{Question}

\begin{document}

\title{On the Generalized Climbing Stairs Problem}

\author{Edray Herber Goins}
\address{Department of Mathematics \\ Purdue University \\ 150 North University Street \\ West Lafayette, IN 47907}
\email{egoins@math.purdue.edu}

\author{Talitha M. Washington}
\address{Department of Mathematics \\ 1800 Lincoln Avenue \\ University of Evansville \\ Evansville, IN 47722}
\email{tw65@evansville.edu}

\begin{abstract}
Let $\mathcal S$ be a subset of the positive integers, and $M$ be a positive integer.   Mohammad K. Azarian, inspired by work of Tony Colledge, considered the number of ways to climb a staircase containing $n$ stairs using ``step-sizes'' $s \in \mathcal S$ and multiplicities at most $M$.

In this exposition, we find a solution via generating functions, i.e., an expression which counts the number of partitions $n = \sum_{s \in \mathcal S} m_s \, s$ satisfying $0 \leq m_s \leq M$.   We then use this result to answer a series of questions posed by Azarian, thereby showing a link with ten sequences listed in the On-Line Encyclopedia of Integer Sequences.  We conclude by posing open questions which seek to count the number of compositions of $n$.
\end{abstract}

\keywords{generating function, partition, sequence, stairs}

\subjclass[2000]{05A15, 05A17}

\maketitle

\section{Introduction}

\begin{quote}
``I'll build a stairway to paradise / with a new step ev'ry day.'' \\ -- Ira and George Gershwin, \emph{Scandals} (1922) \\
\end{quote}

Consider the problem of determining the number of ways to climb a staircase containing $n$ stairs.  For example, there are three ways to climb $n=3$ stairs: take one stair three times, take one stair then two stairs, or take three stairs in one step.  It is natural to restrict this problem by considering the number of ways to climb $n$ stairs where the types of steps taken are either even-numbered or odd-numbered.  For example, there are just two ways to climb a staircase containing $n=3$ stairs if only odd-numbered steps are allowed.  One may further restrict this problem by considering only a distinct set of steps.  For example, there are three such ways to climb a staircase containing $n=3$ stairs: take one stair then two stairs, or take three stairs in just one step.  One may also consider the cases when the order of the steps taken is relevant.

Azarian \cite{Azarian1997A-Generalizatio}, \cite{Azarian2004A-Generalizatio}, inspired by Colledge \cite{Colledge1992Pascals-Triangl}, considered problems which can be phrased as follows:

\begin{staircase} Let $\mathcal S$ be a subset of the positive integers; this will denote the ``sizes'' of steps allowed.  Let $M$ be a positive integer; this will denote the maximum multiplicity of each step-size $s \in \mathcal S$ taken.  What are the possible ways to climb a staircase containing $n$ stairs taking step-sizes $\mathcal S$ and multiplicities at most $M$?
\end{staircase}

\noindent The main focus of this exposition is to give a complete answer to this problem.  In Section 2, we find a solution via generating functions, i.e., an expression which counts the number of partitions $n = \sum_{s \in \mathcal S} m_s \, s$ satisfying $0 \leq m_s \leq M$.   In Section 3, we use this result to answer a series of questions posed by Azarian, thereby showing a link with ten sequences listed in the On-Line Encyclopedia of Integer Sequences.  In Section 4, we conclude by posing open questions which seek to count the number of compositions of $n$.

\section{Main Theorem}

First, we fix some notation.  Let $\mathcal S$ be a subset of the positive integers; this will denote the ``sizes'' of steps taken when climbing the stairs.  For example, $\mathcal S$ can be the set of all even positive integers, the set of all odd positive integers, or perhaps a finite set.  Assume that we have a staircase containing $n$ stairs.  We denote a set of ``steps'' of this staircase by the $\mathcal S$-partition

\begin{equation*} n = \sum_{s \in \mathcal S} m_s \, s \end{equation*}

\noindent where $s$ is the ``step-size'' and $m_s$ denotes the multiplicity.  We will also denote this $\mathcal S$-partition by the sequence

\begin{equation*} \left( m_s \right)_{s \in \mathcal S} \end{equation*}

\noindent where $m_s$ are nonnegative integers, all but finitely many of which are zero.  For example, there are three ways to climb a $n=3$ stairs: take one stair three times, take one stair then two stairs, or take three stairs at once in one step.  If we let $\mathcal S$ denote the set of all positive integers, then we may express these three ways as $\left( 3, 0, 0, 0, \dots \right)$, $\left( 1, 1, 0, 0, \dots \right)$, or $\left( 0, 0, 1, 0, \dots \right)$, respectively.  For a positive integer $M$, denote $p_{\mathcal S}^{(M)}(n)$ as the number of ways to climb a staircase containing $n$ steps using step-sizes $\mathcal S$ and multiplicities at most $M$.  To be precise,

\begin{equation*} p_{\mathcal S}^{(M)}(n)\quad = \quad \# \left \{ \left( m_s \right)_{s \in \mathcal S} \, \left| \, \text{$n = \sum_{s \in \mathcal S} m_s \, s$ \ and \ $0 \leq m_s \leq M$} \right. \right \}. \end{equation*}

\noindent We use the convention $p_{\mathcal S}^{(M)}(0)=1$.

Azarian \cite{Azarian1997A-Generalizatio}, \cite{Azarian2004A-Generalizatio}, inspired by Colledge \cite{Colledge1992Pascals-Triangl}, considered problems which can be phrased as follows:

\begin{staircase} Let $\mathcal S$ be a subset of the positive integers; this will denote the ``sizes'' of steps allowed.  Let $M$ be a positive integer; this will denote the maximum multiplicity of each step-size $s \in \mathcal S$ taken.  What are the possible ways to climb a staircase containing $n$ stairs taking step-sizes $\mathcal S$ and multiplicities at most $M$?
\end{staircase}

\noindent In this case order is irrelevant.  It is easy see that climbing a staircase gives a $\mathcal S$-partition of $n$.   We prove the following as a solution to the problem above:

\begin{partition} Let $p_{\mathcal S}^{(M)}(n)$ denote the number of ways to climb a staircase containing $n$ stairs using step-sizes $\mathcal S$ and multiplicities at most $M$.  Then we have the identity
\begin{equation*} \sum_{n = 0}^{\infty} p_{\mathcal S}^{(M)}(n) \ x^n = \prod_{s \in \mathcal S} \frac {1 - x^{(M+1)s}}{1 - x^s} \quad \text{on the interval} \quad |x| < 1. \end{equation*}
\end{partition}

\noindent \textbf{Remarks.}

\begin{enumerate}
\item[(a)] When $M$ is unbounded (i.e., ``$M = \infty$'') and $\mathcal S$ is the set of all positive integers, then $p_{\mathcal S}^{(M)}(n) = p(n)$ which is the classical partition function.  This result is a generalization of the classical identity found by Euler: \begin{equation*} \sum_{n=0}^{\infty} p(n) \ x^n = \prod_{s=1}^{\infty} \frac 1{1 - x^s} = 1 + x + 2 \, x^2 + 3 \, x^3 + 5 \, x^4 + 7 \, x^5 + \cdots. \end{equation*}
\item[(b)] When $M = 1$ and $\mathcal S$ is the set of all positive integers, then $p_{\mathcal S}^{(M)}(n) = q(n)$ counts the number of ways to climb the staircase using \emph{distinct} step-sizes from $\mathcal S$: \begin{equation*} \sum_{n=0}^{\infty} q(n) \ x^n = \prod_{s=1}^{\infty} \left( 1 + x^s \right) = 1 + x + x^2 + 2 \, x^3 + 2 \, x^4 + 3 \, x^5 + \cdots. \end{equation*}
\end{enumerate}
\noindent These sequences may be verified using the On-Line Encyclopedia of Integer Sequences \cite{oeis}.  Compare these with sequences A000041 and A000009, respectively.

\begin{proof} We review \cite[Theorem 1.1]{MR0557013} and its proof to remind the reader of the language of climbing staircases.  Fix real numbers $r$ and $T$ satisfying $0 < r < 1$ and $0 < T$.  Let $\mathcal S_T \subseteq \mathcal S$ and let $s \in \mathcal S$ such that $s \leq T$.  Define the function

\begin{equation*} G_T(x) = \prod_{s \in \mathcal S_T} \frac {1 - x^{(M+1)s}}{1 - x^s} \quad \text{on the interval} \quad |x| \leq r. \end{equation*}

\noindent Since this function is a finite product of terms involving the geometric series, and the geometric series is uniformly convergent in this compact interval, we may rearrange this expression to yield

\begin{equation*} \begin{aligned} G_T(x) & = \prod_{s \in \mathcal S_T} \frac {1 - x^{(M+1)s}}{1 - x^s} \\ & = \prod_{s \in \mathcal S_T} \left( \frac 1{1 - x^s} - x^{(M+1)s} \cdot \frac 1{1-x^s} \right) \\ & = \prod_{s \in \mathcal S_T} \left( \sum_{m_s = 0}^{\infty} x^{m_s s} - \sum_{m_s=0}^{\infty} x^{(m_s+M+1) s} \right) \\ & = \prod_{s \in \mathcal S_T} \left( \sum_{m_s = 0}^M x^{m_s s} \right). \end{aligned} \end{equation*}

\noindent By comparing the coefficient of $x^n$, we obtain the identity

\begin{equation*} \begin{aligned} G_T(x) & = \prod_{s \in \mathcal S_T} \left( \sum_{m_s = 0}^M x^{m_s s} \right) \\ & = \sum_{n = 0}^{\infty} \# \left \{ \left( m_s \right)_{s \in \mathcal S_T} \, \left| \, \text{$n = \sum_{s \in \mathcal S_T} m_s \, s$ and $0 \leq m_s \leq M$} \right. \right \} x^n. \end{aligned} \end{equation*}

\noindent Since $r$ and $T$ were arbitrary, the theorem holds in the open interval $|x| < 1$ for any given set $\mathcal S$. \end{proof}

\section{Applications}

We are now interested in applying the theorem in the previous section to counting the number of ways to climb a staircase.  In this section, we answer questions of Mohammad Azarian originally posed in \cite{Azarian2004A-Generalizatio}.

\begin{question1}
How many different ways are there to climb a staircase containing $n$ stairs taking distinct even- (odd-, respectively) numbered step-sizes? \end{question1}

Let $M=1$ and $\mathcal S$ be the set of even (odd, respectively) positive integers. Then using the main theorem from the previous section, we have the generating function

\begin{equation*} \sum_{n = 0}^{\infty} p_{\mathcal S}^{(M)}(n) \ x^n = \left \{ \begin{aligned} \displaystyle \prod_{\text{$s$ even}} \left( 1 + x^s \right) & = 1 + x^2 + x^4 + 2 \, x^6 + \cdots \\ &  \text{for even-numbered steps;} \\[10pt] \displaystyle \prod_{\text{$s$ odd}} \left( 1 + x^s \right) & = 1 + x + x^3 + x^4 + x^5 + \cdots \\ &  \text{for odd-numbered steps.} \end{aligned} \right. \end{equation*}

\noindent For the former, compare with sequence A000009 of \cite{oeis}.  For the latter, compare with sequence A000700.

\begin{question2}
For a positive integer $k$, how many different ways are there to climb a staircase containing $n$ stairs taking exactly $k$ stairs at most $k$ times?
\end{question2}

Denote $M = k$ and $\mathcal S = \left \{ k \right \}$.  Then we have the generating function

\begin{equation*} \sum_{n = 0}^{\infty} p_{\mathcal S}^{(M)}(n) \ x^n = \frac {1 - x^{(k+1)k}}{1 - x^k} = 1 + x^k + x^{2k} + \cdots + x^{k^2}. \end{equation*}

\noindent That is, $p_{\mathcal S}^{(M)}(n) \neq 0$ if and only if $n$ is a multiple of $k$ in the form $n = m \, k$ for some positive integer $m \leq k$.

\begin{question3}
How many different ways are there to climb a staircase containing $n$ stairs taking at least two stairs at a time?
\end{question3}

Let $M$ be unbounded and $\mathcal S = \left \{ 2, \, 3, \, 4, \, \dots \right \}$.  Then we have the generating function

\begin{equation*} \sum_{n = 0}^{\infty} p_{\mathcal S}^{(M)}(n) \ x^n = \prod_{s = 2}^{\infty} \frac 1{1 - x^s} = 1 + x^2 + x^3 + 2 \, x^4 + 2 \, x^5 + \cdots. \end{equation*}

\noindent Compare with sequence A002865 of \cite{oeis}.

\begin{question4}
How many different ways are there to climb a staircase containing $n$ stairs taking at most two stairs at a time?
\end{question4}

Let $M$ be unbounded and $\mathcal S = \left \{ 1, \, 2 \right \}$.  Then we have the generating function

\begin{equation*} \begin{aligned} \sum_{n = 0}^{\infty} p_{\mathcal S}^{(M)}(n) \ x^n & = \frac 1{1 - x} \cdot \frac 1{1-x^2} \\ & = 1 + x + 2 \, x^2 + 2 \, x^3 + 3 \, x^4 + 3\, x^5 + \cdots. \end{aligned} \end{equation*}

\noindent Compare with sequence A008619 of \cite{oeis}.

\begin{question5}
How many different ways are there to climb a staircase containing $n$ stairs taking even- (odd-, respectively) numbered step-sizes?
\end{question5}

Let $M$ be unbounded and $\mathcal S$ as the set of even (odd, respectively) positive integers. Then we have the generating function

\begin{equation*} \sum_{n = 0}^{\infty} p_{\mathcal S}^{(M)}(n) \ x^n = \left \{ \begin{aligned} \displaystyle \prod_{\text{$s$ even}} \frac 1{1 - x^s} & = 1 + x^2 + 2 \, x^4 + 3 \, x^6 + \cdots \\ &  \text{for even-numbered steps;} \\[10pt] \displaystyle \prod_{\text{$s$ odd}} \frac 1{1 - x^s} & = 1 + x + x^2 + 2 \, x^3 + 2 \, x^4 + 3 \, x^5 + \cdots \\ &  \text{for odd-numbered steps.} \end{aligned} \right. \end{equation*}

\noindent For the former, compare with sequence A000041 of \cite{oeis}.  For the latter, compare with sequence A000009; it is well-known that $p_{\mathcal S}^{(M)}(n) = q(n)$ in this case.  See \cite{Weis3} for more information.

\begin{question6}
How many different ways are there to climb a staircase containing $n$ stairs where the multiplicity of each step-size is at most 2?
\end{question6}

Denote $M = 2$ and $\mathcal S$ as the  set of positive integers.  Then we have the generating function

\begin{equation*} \begin{aligned} \sum_{n = 0}^{\infty} p_{\mathcal S}^{(M)}(n) \ x^n & = \prod_{s = 1}^{\infty} \left( 1 + x^s + x^{2s} \right) \\ & = 1 + x + 2 \, x^2 + 2 \, x^3 + 4 \, x^4 + 5 \, x^5 + \cdots. \end{aligned} \end{equation*}

\noindent Compare with sequence A000726 of \cite{oeis}.

\begin{question7}
For a positive integer $k$, how many different ways are there to climb a staircase containing $n$ stairs taking exactly $k$ stairs for each step?
\end{question7}

Let $M$ be unbounded and $\mathcal S = \left \{ k \right \}$.  Then we have the generating function

\begin{equation*} \sum_{n = 0}^{\infty} p_{\mathcal S}^{(M)}(n) \ x^n = \frac 1{1 - x^k} = 1 + x^k + x^{2k} + x^{3 k} + x^{4k} + x^{5k} + \cdots. \end{equation*}

\noindent That is, $p_{\mathcal S}^{(M)}(n) \neq 0$ if and only if $n$ is a multiple of $k$.

\begin{question8}
How many different ways are there to climb a staircase containing $n$ stairs where the size of each step is a prime number?  What if the multiplicity of each step-size is at most 1?
\end{question8}

These questions were first considered by Bateman and Erd{\"o}s \cite{MR0079013} and Gupta \cite{MR0078398}, respectively.  Let $M$ be unbounded and $\mathcal S = \left \{ 2, \, 3, \, 5, \, 7, \, \dots, \, \ell, \, \dots \right \}$ denote the set of all prime numbers.  Then we have the generating function

\begin{equation*} \sum_{n = 0}^{\infty} p_{\mathcal S}^{(M)}(n) \ x^n = \prod_{\text{$\ell$ prime}} \frac 1{1 - x^\ell} = 1 + x^2 + x^3 + x^4 + 2 \, x^5 + \cdots. \end{equation*}

\noindent Compare with sequence A000607 of \cite{oeis}.  Now let $M = 1$ with $\mathcal S$ as before.  Then we have the generating function

\begin{equation*} \sum_{n = 0}^{\infty} p_{\mathcal S}^{(M)}(n) \ x^n = \prod_{\text{$\ell$ prime}} \left( 1 + x^\ell \right) = 1 + x^2 + x^3 + 2 \, x^5 + \cdots. \end{equation*}

\noindent Compare with sequence A000586 of \cite{oeis}.

\begin{question9}
How many different ways are there to climb a staircase containing $n$ stairs where the size of each step is a Fibonacci number?  What if the multiplicity of each step-size is at most 1?
\end{question9}

This question was first considered by Klarner \cite{MR0218297}.  Let $M$ be unbounded and $\mathcal S = \left \{ 1, \, 2, \, 3, \, 5, \, \dots, \, F_n, \, \dots \right \}$ denote the set of all Fibonacci numbers, i.e., the Fibonacci sequence with multiplicities removed.  Then we have the generating function

\begin{equation*} \begin{aligned} \sum_{n = 0}^{\infty} p_{\mathcal S}^{(M)}(n) \ x^n & = \prod_{\text{$F_n$ Fibonacci}} \frac 1{1 - x^{F_n}} \\ & = 1 + x + 2 \, x^2 + 3 \, x^3 + 4 \, x^4 + 6 \, x^5 + \cdots. \end{aligned} \end{equation*}

\noindent Compare with sequence A003107 of \cite{oeis}.  Now let $M = 1$ with $\mathcal S$ as before.  Then we have the generating function

\begin{equation*} \begin{aligned} \sum_{n = 0}^{\infty} p_{\mathcal S}^{(M)}(n) \ x^n & = \prod_{\text{$F_n$ Fibonacci}} \left( 1 + x^{F_n} \right) \\ & = 1 + x + x^2 + 2 \, x^3 + x^4 + 2 \, x^5 + \cdots. \end{aligned} \end{equation*}

\noindent Compare with sequence A000119 of \cite{oeis}.

\begin{question10}
For positive integers $a$ and $b$, how many different ways are there to climb a staircase containing $n$ stairs where the size of each step $s$ satisfies $a \leq s \leq b$?
\end{question10}

Let $M$ be unbounded and $\mathcal S$ denote the set integers $s$ satisfying $a \leq s \leq b$.  Then we have the generating function

\begin{equation*} \sum_{n = 0}^{\infty} p_{\mathcal S}^{(M)}(n) \ x^n = \prod_{s=a}^b \frac 1{1 - x^s} = 1 + x^a + \cdots. \end{equation*}

\section{Open Problem}

Now consider the problem of determining the number of ways to climb a staircase containing $n$ stairs -- where the order in which the steps are taken is relevant.  For example, there are \emph{four} ways to climb $n=3$ stairs: take one stair three times, take one stair then two stairs, take two stairs then one stair, or take three stairs at once in one step.  We pose the following problem:

\begin{order} Let $\mathcal S$ be a subset of the positive integers; this will denote the ``sizes'' of steps allowed.  Let $M$ be a positive integer; this will denote the maximum multiplicity of each step-size $s \in \mathcal S$ taken.  What are the possible ways to climb a staircase containing $n$ stairs -- keeping track of the order in which the steps are taken -- if one can only take using step-sizes $\mathcal S$ and multiplicities at most $M$?  
\end{order}

\noindent When the order is irrelevant, climbing a staircase gives a partition of $n$.  However, when the order is relevant, we must consider \emph{compositions} of $n$.

\section*{Acknowledgements}

The authors would like to thank Mohammad K. Azarian for inspiration with this project.


\begin{thebibliography}{1}

\bibitem{MR0557013}
George~E. Andrews.
\newblock {\em The theory of partitions}.
\newblock Addison-Wesley Publishing Co., Reading, Mass.-London-Amsterdam, 1976.
\newblock Encyclopedia of Mathematics and its Applications, Vol. 2.

\bibitem{Azarian1997A-Generalizatio}
Mohammad~K. Azarian.
\newblock A generalization of the climbing stairs problem.
\newblock {\em Mathematics and Computer Education}, 31(1):24--28, Winter 1997.

\bibitem{Azarian2004A-Generalizatio}
Mohammad~K. Azarian.
\newblock A generalization of the climbing stairs problem {I}{I}.
\newblock {\em Missouri Journal of Mathematical Sciences Articles}, 16(1):12--17, Winter 2004.

\bibitem{MR0079013}
P.~T. Bateman and P.~Erd{\"o}s.
\newblock Partitions into primes.
\newblock {\em Publ. Math. Debrecen}, 4:198--200, 1956.

\bibitem{Colledge1992Pascals-Triangl}
Tony Colledge.
\newblock {\em Pascal's Triangle: A Teacher's Guide With Blackline Masters}.
\newblock Tarquin Publications, January 1992.

\bibitem{MR0078398}
Hansraj Gupta.
\newblock Partitions into distinct primes.
\newblock {\em Proc. Nat. Inst. Sci. India. Part. A.}, 21:185--187, 1955.

\bibitem{MR0218297}
David~A. Klarner.
\newblock Representations of {$N$} as a sum of distinct elements from special sequences.
\newblock {\em Fibonacci Quart.}, 4:289--306, 322, 1966.

\bibitem{oeis}
N.~J.~A. Sloan.
\newblock The on-line encyclopedia of integer sequences.
\newblock {\em \newline{\tt http://www.research.att.com/~njas/sequences/}}.

\bibitem{Weis3}
Eric~W. Weisstein.
\newblock Partition function {Q}.
\newblock {\em \newline{\tt http://mathworld.wolfram.com/PartitionFunctionQ.html}}.

\end{thebibliography}
\end{document}